\documentclass[12pt,english,fleqn,liststotoc,bibtotoc,idxtotoc,BCOR7.5mm,tablecaptionabove]{amsart}
\usepackage[T1]{fontenc}
\usepackage[latin1]{inputenc}
\usepackage{geometry}
\usepackage{color}
\usepackage{amsthm}
\usepackage{amstext}
\usepackage{amssymb}
\usepackage{amsmath}
\usepackage{graphicx}
\usepackage{young}
\usepackage[section]{placeins}
\usepackage[svgnames]{xcolor}
\usepackage[unicode=true,
 bookmarks=true,bookmarksnumbered=true,bookmarksopen=false,
 breaklinks=false,pdfborder={0 0 1},backref=false,colorlinks=true]
 {hyperref}
\hypersetup{pdftitle={LyX's Additional Features manual},
 pdfauthor={LyX Team},
 pdfsubject={LyX's additional features documentation},
 pdfkeywords={LyX, Documentation, Additional},
 linkcolor=black, citecolor=black, urlcolor=blue, filecolor=blue,  pdfpagelayout=OneColumn, pdfnewwindow=true,  pdfstartview=XYZ, plainpages=false, pdfpagelabels}
\usepackage{breakurl}

\makeatletter

\theoremstyle{plain}
\newtheorem{thm}{\protect\theoremname}
  \theoremstyle{plain}
  
  \theoremstyle{plain}
  
  \theoremstyle{remark}
  
  \theoremstyle{plain}
  
  \theoremstyle{plain}
  \newtheorem{prop}[thm]{\protect\propositionname}
  \theoremstyle{remark}
  
  \theoremstyle{remark}
  
  \theoremstyle{definition}
  \newtheorem{defn}[thm]{\protect\definitionname}
  \theoremstyle{definition}
  \newtheorem{example}[thm]{\protect\examplename}
  \theoremstyle{plain}

\usepackage{ifpdf}
\ifpdf

 \IfFileExists{lmodern.sty}
  {\usepackage{lmodern}}{}

\fi 

\usepackage{multicol}

\makeatother

  \providecommand{\claimname}{\inputencoding{latin9}Claim}
  \providecommand{\conjecturename}{\inputencoding{latin9}Conjecture}
  \providecommand{\corollaryname}{\inputencoding{latin9}Corollary}
  \providecommand{\definitionname}{\inputencoding{latin9}Definition}
  \providecommand{\examplename}{\inputencoding{latin9}Example}
  \providecommand{\lemmaname}{\inputencoding{latin9}Lemma}
  \providecommand{\notename}{\inputencoding{latin9}Note}
  \providecommand{\propositionname}{\inputencoding{latin9}Proposition}
  \providecommand{\questionname}{\inputencoding{latin9}Question}
  \providecommand{\remarkname}{\inputencoding{latin9}Remark}
\providecommand{\theoremname}{\inputencoding{latin9}Theorem}

\newcommand{\red}{\textcolor{red}}

\begin{document}
\author{Nick Early}
\thanks{Department of Mathematics, Penn State University, University Park, PA 16802; \\
	earlnick@gmail.com}

\title[Representation Theory for Generalized Permutohedra: Simplicial Plates]{Combinatorics and Representation Theory for Generalized Permutohedra I: Simplicial Plates}
\maketitle
\begin{abstract}

In this paper, we announce results from our thesis, which studies for the first time the categorification of the theory of generalized permutohedra.  The vector spaces in the categorification are tightly constrained by certain continuity relations which appeared in physics in the mid 20th century.  We describe here the action of the symmetric group on the vector spaces in this categorification.  Generalized permutohedra are replaced by vector spaces of characteristic functions of polyhedral cones about faces of permutohedra, called \textit{plates}, due to A. Ocneanu.  The symmetric group acts on plates by coordinate permutation.  

In combinatorics, the Eulerian numbers count the number of permutations with given numbers of ascent and descents.  The classical Worpitzky identity expands a power $r^p$ as a sum of Eulerian numbers, with binomial coefficients.  In our thesis, for the main result we generalize the classical Worpitzky identity to an isomorphism of symmetric group modules, corresponding geometrically to the tiling of a scaled simplex by unit hypersimplices.  In the categorification, the volume of a hypersimplex is replaced by the complex-linear dimension of a vector space associated to it. The main technical aspect of the proof of the character formula for the simplex involves a partition of unity of a commutative algebra of translations on a discrete torus, and a certain modular Diophantine equation.

A detailed paper is in preparation.

\end{abstract}

\begingroup
\let\cleardoublepage\relax
\let\clearpage\relax
\tableofcontents
\endgroup

\newpage

\section{Introduction}

This paper is an extension of a program on higher representation theory originated by Ocneanu a long time ago.  At the core of his model is the so-called all-subset hyperplane arrangement, consisting of all hyperplanes $\{x\in\mathbb{R}^n\sum_{i\in I}x_i,\ \sum_{i=1}^n x_i=0\}$ for $I$ a proper subset of $\{1,\ldots, n\}$.  The all-subset hyperplane arrangement is actually known in various areas, but perhaps due to its complexity and the lack of a systematic approach, it has not directly received a concentrated treatment in mathematics.  The all-subset hyperplane arrangement dates back at least to the 1960's in works of H. Araki and D. Ruelle on Wightman functions and generalized retarded Green functions in Quantum Field Theory.  See for example \cite{Streater} for an overview.

We study the action of the symmetric group on certain finite-dimensional vector spaces constructed from families of convex polyhedral regions in Euclidean space, called \textbf{plates}.  The vector spaces are the linear spans of the characteristic functions of families of plates.

Our main result is a proof of a formula for the characters of the symmetric group that we obtain in this way.  It generalizes a well-known formula in combinatorics, called the Worpitzky identity, that from our point of view is obtained combinatorially by evaluating characters on the identity permutation.  It corresponds geometrically to the decomposition of a scaled simplex into unit hypersimplices.

The paper is structured as follows.  We shall begin by describing the Worpitzky identity and relevant surrounding ideas in combinatorics.  Then we shall summarize results from \cite{OcneanuPlates} on the essential properties of plates, along with several examples and illustrations.  Following that we formulate our main result, the character formula, along with some further results.  

\subsection{Background}

In \cite{PostnikovPermutohedra}, A. Postnikov introduced and studied generalized permutohedra.  By \cite{OcneanuPlates}, as plates are characteristic functions of cones around faces of permutohedra, it follows from an inclusion/exclusion argument that any generalized permutohedron can be realized from a linear combination of plates with coefficients $\pm1$.  Replacing generalized permutohedra with vector spaces of characteristic functions of cones amounts to a categorification.  We study the combinatorics and representation theory of this categorification.

Let us point out some subjects in Mathematics which are equivalent or related to the theory of plates, along with suggested references.

\begin{enumerate}
	\item Matroids and polymatroids.
	\begin{enumerate}
		\item Plates in hypersimplices $B_{a,b}$ are closely related to matroids, as in \cite{GelfandMacPherson,MatroidsGrassmannian}.  In \cite{OcneanuPlates} Ocneanu develops elaborate combinatorial techniques which show how to localize generalized permutohedra to hypersimplices.  It is shown that these localizations are matroid polytopes, and any matroid polytope is a linear combination of such localizations.
		\item Plates in simplices $\Delta_r^n$ are closely related to polymatroids, in the sense of \cite{DerksenFink}.
	\end{enumerate}

	\item The module of nondegenerate plates about a point is isomorphic to the free Lie module.
	\begin{enumerate}
		\item In classical representation theory, according to the standard construction one works in the regular representation to obtain the irreducible representations of the symmetric group, by projecting with the Young idempotents onto the irreducible submodules.  It is well-known that the regular representation of $S_{n-1}$ is the restriction of an $S_n$-module called the \textbf{free Lie module}, see \cite{ReutenauerFreeLieAlgebras}.  It is spanned by the $n!$ alternating binary trees, or bracketings, where each leaf is labeled with an integer $1,\ldots, n$, and has basis the $(n-1)!$ ``combed'' antisymmetric trees $\lbrack \ldots\lbrack1,i_2,\rbrack, i_3\rbrack,\ldots\rbrack, i_{n-1}\rbrack ,i_n\rbrack$.
		\item The fundamental observation in \cite{OcneanuPlates} is that the module of plates about a point is isomorphic to the free Lie module in brackets of $n$ (distinct) generators, modulo certain relations coming from plates with fewer lumps.  Here $S_n$ acts on the space of complex-linear combinations of plates by permuting variables.
		
		Here the $S_n$-module of nondegenerate plates has basis the $(n-1)!$ permutations of $2,3,\ldots, n$ of the standard plate 
		$$\lbrack\lbrack 1_0 2_0\cdots n_0\rbrack\rbrack,$$
		which has support the cone generated by the simple $A_{n-1}$ roots $e_i-e_{i+1}$ extending out from the point $(0,\ldots, 0)$.  Note that the $S_n$-action crucially fixes the position of the plate, $(0,\ldots, 0)$, which would not be the case for a plate with position $(s_1,\ldots, s_n)$, where not all $s_i$ are the same.
	\end{enumerate}
	\item Descent statistics for permutations and toric geometry.
	\begin{enumerate}
		\item The $S_n$ hypersimplex plate modules $\text{Pl}\left(B_{a,b}\right)$ are isomorphic to the modules in \cite{HendersonWachs1} which describe the $S_n$-action on the cohomology of a toric variety associated to a Coxeter arrangement \cite{StembridgeToric}.  In particular, the character values for plates in hypersimplices can be obtained as coefficients of the quasi-Eulerian symmetric functions which appear in \cite{HendersonWachs1}.  See \cite{DolgachevLunts,ProcesiToric,StembridgeToric} for the related toric geometry and \cite{Stapledon} for more recent developments.
	\end{enumerate}
	\item Word-quasi-symmetric functions and \'Ecalle's mould calculus.
	\begin{enumerate}
		\item In \cite{ThibonMould} rational functions in $n$ formal indeterminants $y_1,\ldots, y_n$ are given to provide a functional model for the so-called word-quasisymmetric functions, WQSym.  These functions, which are given as examples of \'Ecalle's so-called moulds, are closely related to plates.  Details will appear elsewhere.

	\end{enumerate}

\end{enumerate}

\subsection{The Classical Worpitzky Identity}

The Eulerian numbers $E_{i,j}$, or $A_{n,j}$ in the usual notation, were defined by Euler in 1755.  They count permutations of $\{1,\ldots, n\}$ with $i$ ascents and $j=(n-1)-i$ descents.  Our notation differs slightly from the standard one in order to emphasize a correspondence with plates in hypersimplices.  The following table starts at the top with $E_{0,0}$ and continues in the second row with $E_{0,1},E_{1,0},\ldots$.
$$\begin{array}{ccccccccc}
&  &  &  & 1 &  &  &  &  \\ 
&  &  & 1 &  & 1 &  &  &  \\ 
&  & 1 &  & 4 &  & 1 &  &  \\ 
& 1 &  & 11 &  & 11 &  & 1 & \\ 
1 &  & 26 &  & 66 &  & 26 &  & 1\\
&&&&\vdots &&&&
\end{array}$$
For example, the third row counts the permutations
$$321,\ \ 132,213,231,312,\ \ 123.$$
In 1883 Worpitzky discovered an identity which expresses a power as a sum of Eulerian numbers with binomial coefficients.  For example,
\begin{eqnarray*}
	& x^{2}=\binom{x}{2}\cdot 1+\binom{x+1}{2}\cdot 1&\\
	& x^{3}=\binom{x}{3}\cdot 1+\binom{x+1}{3}\cdot 4+\binom{x+2}{3}\cdot 1 &\\
	& \vdots &
\end{eqnarray*}
See \cite{Petersen} for a detailed discussion of the Eulerian numbers and many related topics.

We generalize this combinatorial identity to an identity of characters of the symmetric group $S_n$.  These are characters of plate modules for the simplex and hypersimplices in $n$ coordinates.  A detailed paper is in preparation.

\section{Plates}
From this section on, by a \textbf{plate} we shall mean the \textbf{characteristic function} of the polyhedral cone defined by a flag of inequalities of the form

\begin{eqnarray*}
	x_{S_1} & \ge & s_1\\
	x_{S_1}+x_{S_2} & \ge & s_1+s_2\\
	& \vdots & \\
	x_{S_1} + \cdots+x_{S_{k-1}}& \ge & s_1+\cdots +s_{k-1}\\
	x_{S_1} + \cdots+x_{S_{k}}& = & s_1+\cdots +s_{k},
\end{eqnarray*}
where the last line is an equality.  Here $x_S=\sum_{i\in S}x_i$, where $S$ is a subset of $\{1,\ldots, n\},$ with $x_1,\ldots, x_n$ real variables and  $S_1\sqcup\cdots\sqcup S_n=\{1,\ldots, n\}$ is an ordered set partition and $s_1+\cdots+s_k=r\in\mathbb{N}$ is an ordered partition.  In most of what follows we will assume $x_i\ge0$ and $r\ge 1$.  The sets $S_i$ are called \textbf{lumps} $S_i$ and the corresponding integers $s_i$ are called \textbf{positions}.  We will use the notation $\lbrack\lbrack(S_1)_{s_1}\cdots (S_k)_{s_k}\rbrack\rbrack$ to indicate the characteristic function of the region determined by the above equations.
\begin{figure}[h!]
	\centering
	\includegraphics[width=0.75\linewidth]{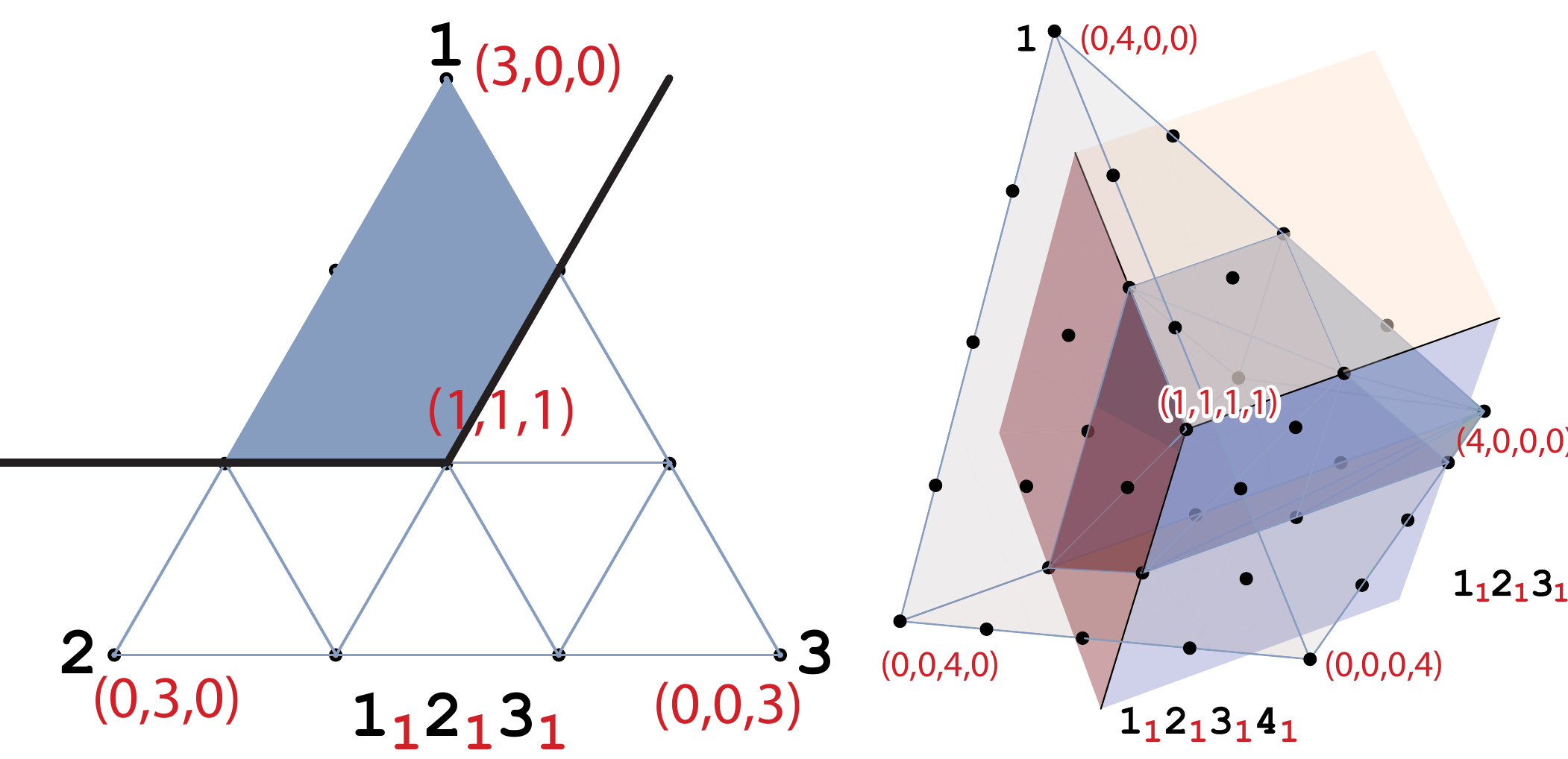}
	\caption{Standard nondegenerate plates as polyhedral cones, generated by simple roots.}
	\label{fig:standardplates}
\end{figure}
In more detail, a plate is the characteristic function of the cone about a face of the permutohedron, the polytope formed by the convex hull of permutations of $(0,1,2,\ldots, n-1)\in\mathbb{R}^n$.  In Figure \ref{fig:standardplates} the standard nondegenerate plates, i.e. all lumps have size 1, are given for $n=3$ and $4$ as unbounded polyhedral cones generated by simple roots of type $A_2$ respectively $A_3$.  For example, in the $n=4$ coordinate case on the right, the plate is identified with the region
$$\{(1,1,1,1)+t_1\cdot (1,-1,0,0)+t_2(0,1,-1,0)+t_3(0,0,1,-1):t_i\ge 0\},$$
opening toward the vertex $(4,0,0,0)$.

The permutohedron itself was studied and generalized previously in \cite{PostnikovPermutohedra}.  For nondegenerate plates, as in for example Figure \ref{fig:standardplates}, from the inequalities above one can verify that the plate is the convex span of a system of simple roots of type $A_{n-1}$.  

The nondegenerate plate $\lbrack\lbrack1,2,\ldots, n\rbrack\rbrack$ was essentially studied before in \cite{GelfandGraevPostnikov}, where a finite part of it is constructed as the convex hull of the positive roots $e_i-e_j$ of type $A_{n-1}$, for $1\le i<j\le n$.

\subsection{Plate Relations}

In \cite{OcneanuPlates}, Ocneanu establishes fundamental properties of plates.  We shall need a subset of these results, which we reproduce without proof in Theorems \ref{CyclicSum}, \ref{StandardBasis} and \ref{PlateRelations}.

\begin{defn}
	Let $\pi,\pi'$ be plates.  We shall say that $\pi'$ is a lumping of $\pi$ if for each lump $A'_{a'}$ of $\pi'$, $A'$ is a union of consecutive lumps of $\pi$ and if $a'$ is the sum of their respective positions.
\end{defn}

For example, $\pi'=\lbrack\lbrack 1_a 24_{b+c}3_d\rbrack\rbrack$ is a lumping of $\pi=\lbrack\lbrack 1_a 2_b 4_{c}3_d\rbrack\rbrack$.  

The lumping of $x_2$ and $x_4$ into $x_{24}=x_2+x_4$ corresponds to a projection of the ambient simplex along the edge with equation $x_2+x_4=b+c$, onto a lower dimensional space parametrized by the variables $x_1,x_{24},x_3$.  Under this projection, the two simplex vertices labeled with $2$ and $4$, which are the endpoints of the projection edge, are identified.

As plates in $n$ variables are characteristic functions of flags of at most $n-1$ inequalities, intersected with a hyperplane $\sum x_i=r$, we see now that lumping can be understood geometrically directly in terms of the defining equations.

\begin{prop}
	The lumpings of a plate $\pi$ are obtained from $\pi$ by removing one or more bounding hyperplanes, i.e. we delete a subset of the $\le n-1$ inequalities which define $\pi$.
\end{prop}

According to Theorem \ref{CyclicSum}, a plate $\pi$ with $l$ lumps determines a decomposition of its ambient space into a sum of the $l$ cyclic rotations of $\pi$ which intersect only on their faces.
\begin{thm}\label{CyclicSum}
	Let $r=a+b+\cdots+c$ and $n\ge 1$ be given.  Let $(A,B,\ldots, C,D)$ be an ordered set partition of $\{1,\ldots, n\}$.  Then we have the cyclic sum relation 
	$$\lbrack\lbrack ({S_1S_2\cdots S_{k-1}S_k})_{r}\rbrack\rbrack=\lbrack\lbrack (S_1)_{s_1} (S_2)_{s_2}\cdots (S_{k-1})_{s_{k-1}} (S_k)_{s_{k}}\rbrack\rbrack+\lbrack\lbrack  (S_{k})_{s_{k}}(S_1)_{s_1} (S_2)_{s_2}\cdots (S_{k-1})_{s_{k-1}}\rbrack\rbrack+\cdots$$
	$$+ \lbrack\lbrack  (S_{2})_{s_{2}} \cdots (S_{k-1})_{s_{k-1}}(S_k)_{s_k}(S_{1})_{s_{1}}\rbrack\rbrack.$$
\end{thm}

In practice, all computations involving plates are made possible by the formula given in Theorem \ref{PlateRelations} below.  In particular, Theorem \ref{StandardBasis} provides a basis, and using Theorem \ref{PlateRelations} one can build matrices for the action of permutations.  To prove our character formula, however, we introduce a basis which is permutation invariant up to a root of unity.  Using this basis, we parametrize explicitly the diagonal of the matrix of any given permutation, from which it is immediate to obtain its trace.  
This basis, the \textbf{q-basis}, inspires an interesting algebraic structure which we introduce in Definition \ref{PermutohedronAlgebra}.

We come now to some of the main results of \cite{OcneanuPlates}, stated in what follows without proof.

\begin{thm}\label{StandardBasis}
	The set of all plates which have $1$ in the first lump is linearly independent.  In particular, it is a basis.
\end{thm}
In what follows, we call the basis of Theorem \ref{StandardBasis} the \textbf{standard basis}.

In \cite{OcneanuPlates}, Ocneanu proves the fundamental Theorem \ref{PlateRelations}, giving the module of relations for plates, using homological arguments involving rooted trees.  

The proof involves properties of hypergeometric functions of Euler Beta type for special parameter configurations.

\begin{thm}\label{PlateRelations}
	Let
	$$\pi =\lbrack\lbrack(S_m)_{s_m}(S_{m-1})_{s_{m-1}}\cdots (S_1)_{s_1}(S_{m+1})_{s_{m+1}}\cdots (S_k)_{s_k}\rbrack\rbrack,$$
	labeled so that $1\in S_1$.  This expands in the standard basis as
	$$\pi=\sum_{\pi'\in\text{shL}\left(((S_1)_{s_1},\ldots, (S_m)_{s_m}),((S_{m+1})_{s_{m+1}},\ldots, (S_{k})_{s_k})\right)}(-1)^{m-1}(-1)^{k-n_L}\lbrack\lbrack\pi'\rbrack\rbrack,$$
	where $n_L$ is the number of lumps of $\pi'$ and where the sum is over all lumped shuffles 
	$$(S_{1},S_{i_2},\ldots, S_{i_k})$$
	of $(S_1,\ldots, S_m)$ and $ (S_{m+1},\ldots, S_{k})$ so that the lumps $S_{m+1},\ldots, S_{k}$, i.e. those after $S_1\ni 1$, are not lumped together.
\end{thm}

\begin{example}
	We give the lumped shuffle expansions of two plates in the standard basis.
	$$\lbrack\lbrack 35_a 124_b 6_c\rbrack\rbrack=\lbrack\lbrack124_b 356_{a+c}\rbrack\rbrack + \lbrack\lbrack12345_{a+b} 6_c\rbrack\rbrack - 
	\lbrack\lbrack124_b 6_c 35_b\rbrack\rbrack - \lbrack\lbrack124_a 35_b 6_c\rbrack\rbrack$$
\end{example}
$$\lbrack\lbrack35_b 124_a 6_c 7_d\rbrack\rbrack=\lbrack\lbrack 124_a6_c 3 57_{b+d}\rbrack \rbrack + \lbrack\lbrack 124_{a}356_{a+c}7_d\rbrack \rbrack + \lbrack\lbrack 12345_{a+b}, 6_{c} 7_d\rbrack \rbrack$$
$$ - \lbrack\lbrack 124_{a}6_{c} 7_d 35_b\rbrack \rbrack - \lbrack\lbrack 124_{a}6_c 35_b 7_d\rbrack \rbrack - \lbrack\lbrack 124_{a} 35_b 6_c 7_d\rbrack \rbrack$$

\begin{example}
	In Figure \ref{fig:platerelations} two more plate relations are given.
	\begin{figure}[h!]
		\centering
		\includegraphics[width=1\linewidth]{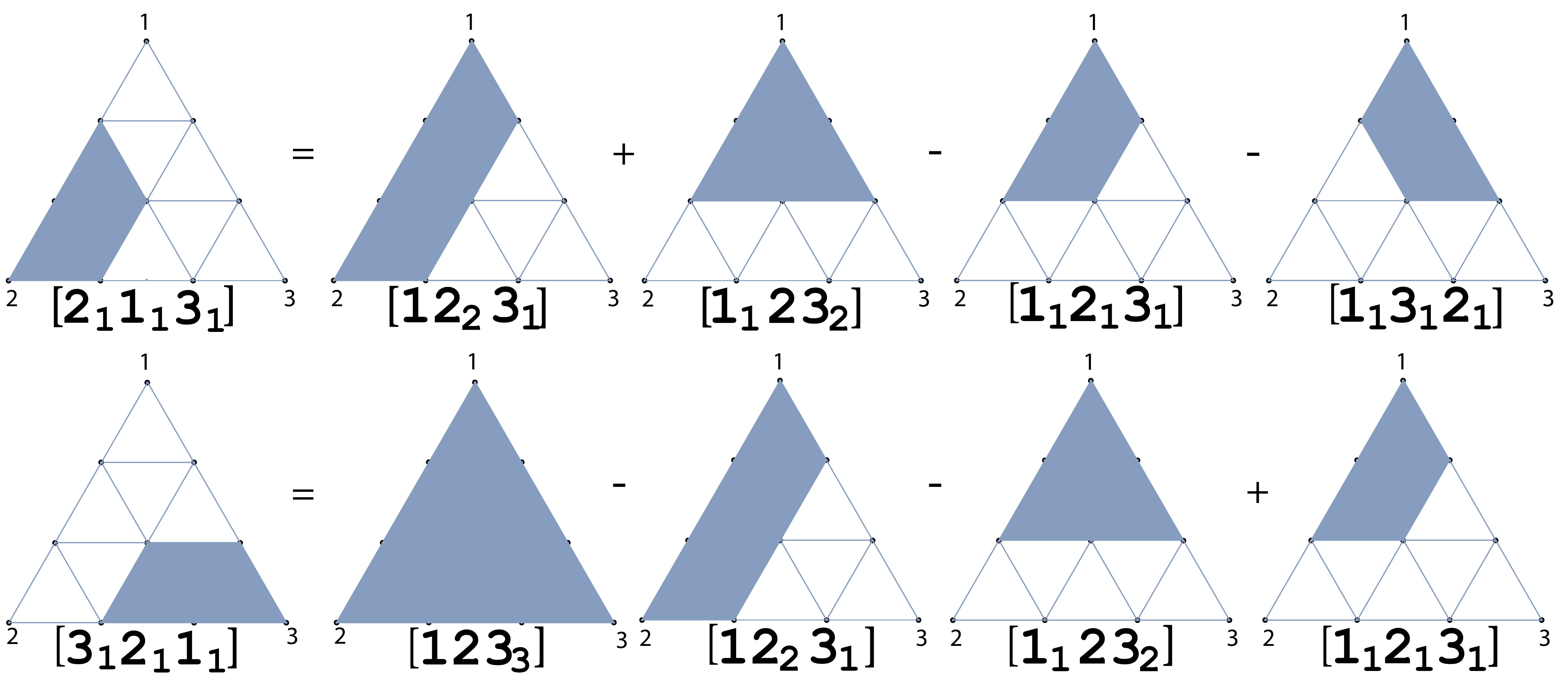}
		\caption{}
		\label{fig:platerelations}
	\end{figure}
\end{example}

\section{Worpitzky Isomorphism}

For positive integers $a,b$ with $a+b=n$, let
$$B_{a,b}=\{x\in\lbrack0,1\rbrack^{n}:\sum_{i=1}^{n} x_i=a\},$$
and for each $r\ge 1$ let
$$\Delta_r^n=\{x\in\lbrack0,r\rbrack^n:\sum_{i=1}^n x_i=r\}.$$
Here $B_{a,b}$ is the well-known hypersimplex, equivalently realized as the convex hull of permutations of the vector $(1,\ldots, 1,0,\ldots, 0)$ with $a$ 1's and $b$ 0's.

We denote by $\text{Pl}(\Delta_r^n)$ and $\text{Pl}(B_{a,b})$ the complex-linear spans of plates which have support in $\Delta_r^n$ and respectively $B_{a,b}$.

\begin{figure}[h!]
	\centering
	\includegraphics[width=0.6\linewidth]{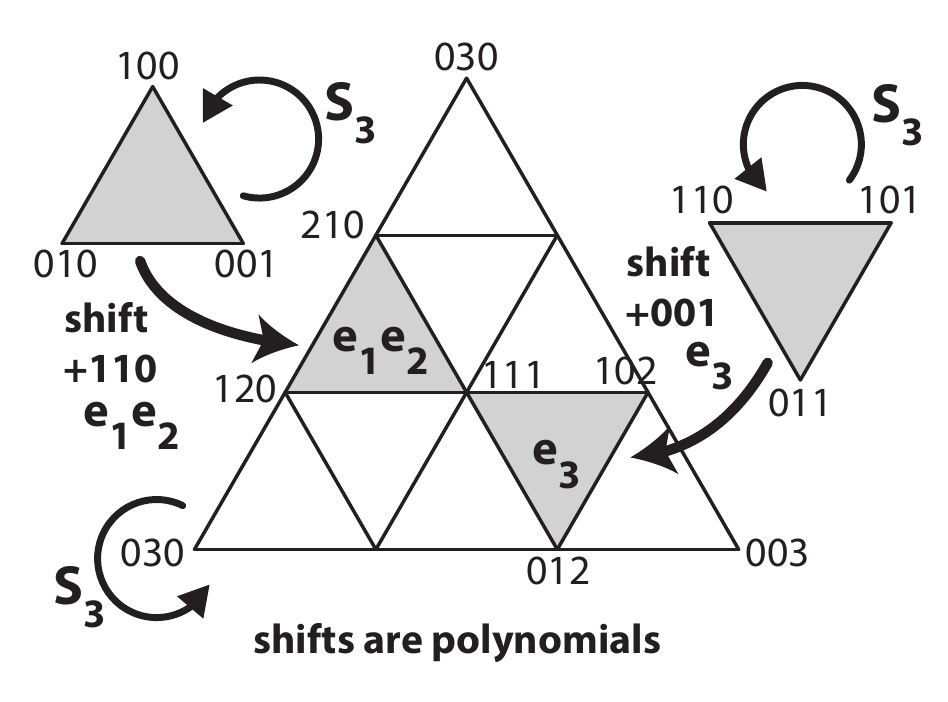}
	\caption{Polynomials Encode Position}
	\label{fig:PositionAsPolynomial}
\end{figure}

\begin{example}
	Figure \ref{fig:PositionAsPolynomial} illustrates the $n=3$ case for $\text{Pl}(\Delta_3^3)$.  Here, in $n=3$ variables, there are two kinds of hypersimplices, namely the ``down triangle'' $B_{1,2}$ and the ``up triangle'' $B_{2,1}$.  Ocneanu shows constructively in \cite{OcneanuPlates} that, in general, $\text{Pl}(B_{a,b})$ is a submodule of $\text{Pl}(\Delta_a^{a+b})$
	
	The hypersimplex plate $\lbrack\lbrack123_1\rbrack\rbrack\in\text{Pl}(B_{1,2})$ is the convex hull of $(1,0,0),(0,1,0),(0,0,1)$, and $\lbrack\lbrack123_2\rbrack\rbrack\in\text{Pl}(B_{2,1})$ is the convex hull of $\{(1,1,0),(0,1,1),(1,0,1)\}$.  The plates $\lbrack\lbrack123_1\rbrack\rbrack$ and $\lbrack\lbrack123_2\rbrack\rbrack$ are thus identified with $B_{1,2}$ and respectively $B_{2,1}$.
	
	In Figure \ref{fig:PositionAsPolynomial}, the hypersimplex plates $\lbrack\lbrack123_1\rbrack\rbrack$ and $\lbrack\lbrack123_2\rbrack\rbrack$
	are translated by the vectors $(1,1,0)$ and respectively $(0,0,1)$.  The position data can be conveniently stored using the monomials, or symmetric tensors, $(e_1 e_2)$ and $e_3$, where $\{e_1,e_2,e_3\}$ is the standard basis for $\mathbb{C}^3$.  
\end{example}

\begin{thm}\label{EquivariantWorpitzky}
There is an isomorphism of $S_n$ modules $$\text{Pl}(\Delta_r^{n})\simeq \bigoplus_{a=1}^{n-1} \text{Sym}^{r-a}(\mathbb{C}^n)\otimes \text{Pl}(B_{a,n-a}),$$
where by convention $\text{Sym}^k(\mathbb{C}^n)=0$ if $k<0$.
\end{thm}
This follows from a straightforward geometric argument together with a result of Ocneanu \cite{OcneanuPlates} that 
$$\dim(\text{Pl}(B_{a,b}))=E_{a-1,n-a-1}.$$
The characters of the identity permutation (hence the dimensions of $\text{Pl}(\Delta_r^{n}),$ $\text{Sym}^{r-a}(\mathbb{C}^n)$, and $ \text{Pl}(B_{a,n-a})$) satisfy the classical Worpitzky identity
$$r^{n-1}=\sum_{a=1}^{n-1}\dim(\text{Sym}^{r-a}(\mathbb{C}^n))\dim(\text{Pl}(B_{a,n-a}))= \sum_{a=1}^{n-1} {n+r-a-1\choose n-1}E_{a-1,n-a-1}.$$

We have thus \textbf{categorified} the classical Worpitzky identity.  We replace a positive counting formula with an identity of characters of symmetric group representations, realized concretely as complex linear spaces spanned by characteristic functions of regions called plates.  Our vector spaces, however, are tightly constrained.  They capture the ``linear'' properties of quantum field theory, see for example \cite{Streater}.  We shall discuss this further in a subsequent paper.

\section{More on Plates}

Recall that a \textbf{basis plate} is a plate $\lbrack\lbrack (S_1)_{s_1} (S_2)_{s_2}\cdots (S_{k-1})_{s_{k-1}} (S_k)_{s_k}\rbrack\rbrack$ such that $1\in S_1$, which means geometrically that it contains the direction towards the vertex labeled 1.  A quantum or \textbf{q-plate}, which assigns a root of unity to each point in space, has the notation ``$\{\cdots \}$'' and is defined by the cyclic sum
\begin{eqnarray*}
	\{(S_1)_{s_1}\cdots (S_k)_{s_k}\} & = & q^0\lbrack\lbrack (S_1)_{s_1} (S_2)_{s_2}\cdots (S_{k-1})_{s_{k-1}} (S_k)_{s_k}\rbrack\rbrack+q^{-s_k}\lbrack\lbrack (S_k)_{s_k} (S_1)_{s_1} (S_2)_{s_2}\cdots C_c\rbrack\rbrack\\
	& + & q^{-(s_{k-1}+s_{k})}\lbrack\lbrack (S_{k-1})_{s_{k-1}}(S_{k})_{s_{k}} (S_{1})_{s_{1}}(S_{2})_{s_{2}} \cdots \rbrack\rbrack +\cdots
\end{eqnarray*}
in which the power of $q$ is the $(-1)$ times the sum of the positions of the lumps moved to the front.

A \textbf{basis q-plate} is a q-plate with the sum normalized to have coefficient $q^0=1$ when $1$ is in the first lump, in which case we have the basis q-plate.  Note that q-plates are basis q-plates, up to scaling by a root of unity.

Examples are given below to illustrate what is going on geometrically.

\begin{figure}[h!]
	\centering
	\includegraphics[width=.9\linewidth]{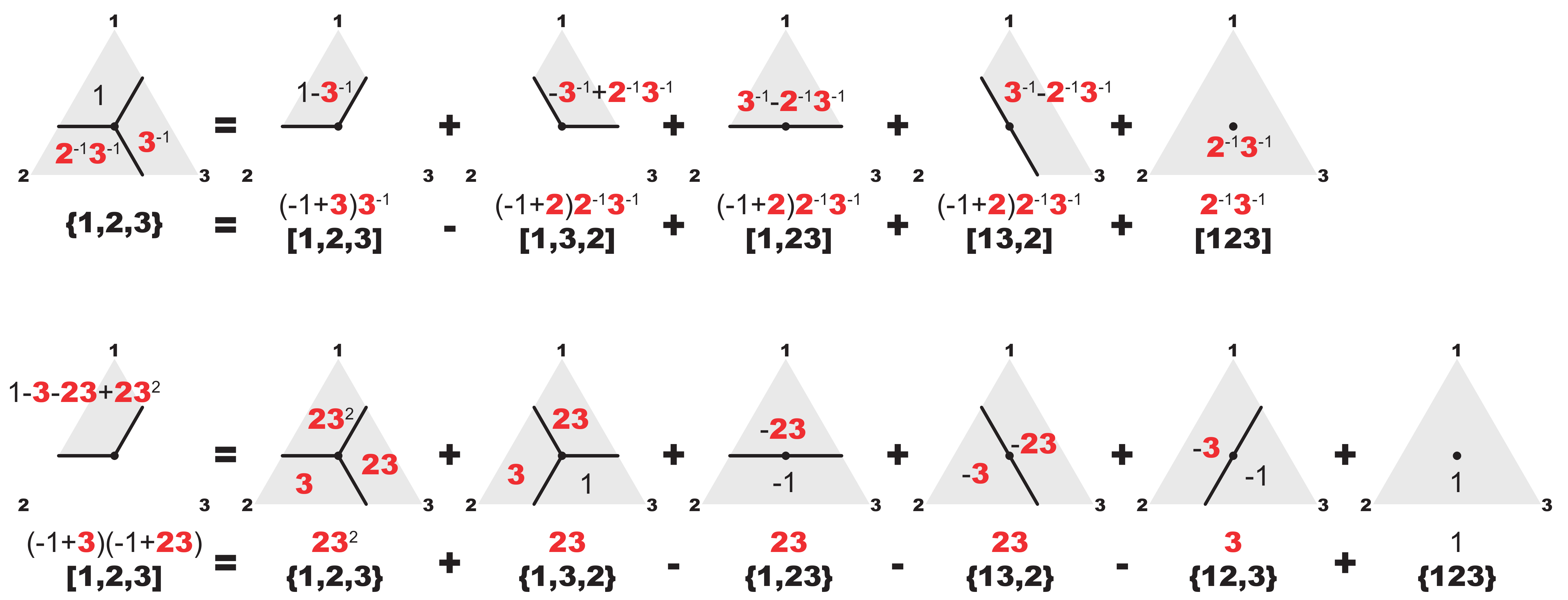}
	\caption{Basis change: q-Plate $\Leftrightarrow$ plate}
	\label{fig:qplatechangebasis0}
\end{figure}
\begin{figure}[h!]
	\centering
	\includegraphics[width=0.45\linewidth]{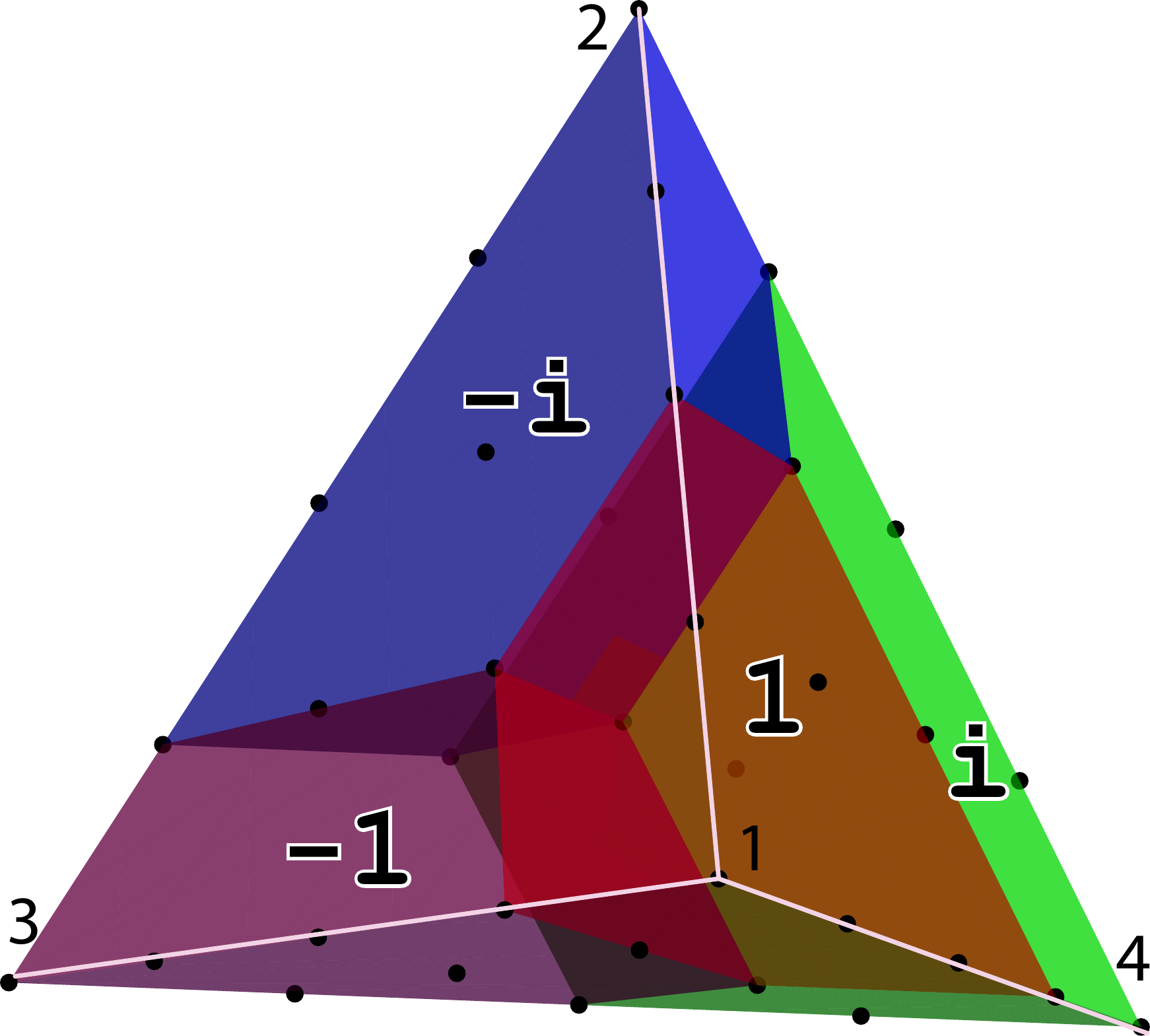}
	\caption{q-plate $\{1_1 2_1 3_1 4_1\}$ in dimension 3}
	\label{fig:qplate3d}
\end{figure}

In Figure \ref{fig:qplatechangebasis0} we use the notation $\mathbf{\red{1}}=q^{a}$, $\mathbf{\red{2}}=q^{b}$ and  $\mathbf{\red{3}}=q^{c}$, where $q=e^{-2\pi i/r}$ and $r=a+b+c$, to describe both the coefficients in the expansion of the q-plate $\{1_a 2_b 3_c\}$ in terms of basis plates and conversely.  Note that $\mathbf{\red{1}}$ does not appear due to the relation $ q^a q^b q^c=1$.  The subscripts $a,b,c$ are implicit in Figure \ref{fig:qplatechangebasis0}.  

Figure \ref{fig:qplate3d} depicts the graph of the q-plate 
$$\{1_1 2_1 3_1 4_1\}=\lbrack\lbrack1_1 2_1 3_1 4_1\rbrack\rbrack+q^{-1}\lbrack\lbrack4_1 1_1 2_1 3_1\rbrack\rbrack+q^{-(1+1)}\lbrack\lbrack3_1 4_11_1 2_1\rbrack\rbrack+q^{-(1+1+1)}\lbrack\lbrack2_1 3_1 4_1 1_1\rbrack\rbrack$$
as a cyclic sum of four plates in the simplex $\Delta_4^4$, each weighted with a power of $q=e^{-2\pi i/4}$.

\section{Formulation of the Main Theorem}

There was a conjecture due to Ocneanu about the characters of the $S_n$-representation of plates as functions on the simplex $\Delta_r^n$.  The characters are remarkable due to the fact that the values are natural numbers with number-theoretic properties.  The conjecture, which we prove in our thesis, is the following.

\begin{thm}
	If $\sigma$ is a permutation with cycle lengths $\lambda_1,\ldots, \lambda_k$, acting on plates in $\Delta_r^n$, then the character value of $\sigma$ is $r^{k-1}$ if $\gcd(\lambda_1,\ldots, \lambda_k,r)=1$ and 0 otherwise.
\end{thm}
\subsection{Method of Proof}

The logical structure of the proof, which will appear separately, is presented in the diagram below, where each $\Updownarrow$ represents an isomorphism.  The aim is to prove the character formula for the bottom row, for plates.  The proof involves counting solutions to a Diophantine equation in the top row.

The bottom double arrow $\Updownarrow$ involves a quite intricate algorithm due to Ocneanu and will appear in our joint work, \cite{EarlyOcneanu}.

\begin{eqnarray*}
	& \mathcal{I}_r^n=\{I\in\left(\mathbb{Z}\slash r\right)^n:i_1+\cdots+i_n\equiv 1\text{ mod r}\}&\\
	& \Updownarrow&\\
	&\text{Partition of unity of } \mathcal{C}_r^n\text{, with idempotents } \epsilon_I=\frac{1}{r^n}\prod_{k=1}^n\left(\sum_{j=0}^{r-1}e^{\left(-2\pi \sqrt{-1}ji_k/r\right)}e_k^j\right)&\\
	&\Updownarrow&\\
	&\text{Monomial basis of } \mathcal{C}_r^n:  \{e_2^{l_2}\cdots e_n^{l_n}:l_i=0,\ldots, r-1\}&\\
	&\Updownarrow&\\
	&\text{q-plates}: \{(S_1)_{s_1}\cdots (S_k)_{s_k}\}&\\
	&\Updownarrow&\\
	&\text{Plates}: \lbrack\lbrack (T_1)_{t_1}\cdots (T_l)_{t_l}\rbrack\rbrack&
\end{eqnarray*}

\subsection{Algebra of Translations}
We introduce the basic terminology and statements of results, postponing the details to a future publication.
\begin{defn} \label{PermutohedronAlgebra}
Let $q=e^{-2\pi i/r}$.
	The translation algebra $\mathcal{C}_r^n$ is the commutative algebra over $\mathbb{C}$ given in terms of generators and relations as
	$$\mathcal{C}_r^n=\left\langle e_1,\ldots, e_n:e_i^r=1, e_i e_j=e_j e_i,e_1\cdots e_n=q \right\rangle.$$
\end{defn}

\begin{prop}
	The symmetric group $S_n$ acts on $\mathcal{C}_r^n$ by permuting the variables.  Moreover, as a complex vector space, $\mathcal{C}_r^n$ has dimension $r^{n-1}$, and the set 
	$$\{e_2^{j_2}\cdots e_n^{j_n}:0\le j_i\le r-1\}$$ 
	consisting of monomials normalized so that $e_1$ does not appear, is a basis.
\end{prop}
A remarkable property of $\mathcal{C}_r^n$ is that it acts freely and transitively on $\text{Pl}(\Delta_r^n)$.  

\begin{thm}\label{permutohedronAlgIsomorphism}
	The algebra $\mathcal{C}_r^n$ and the simplex plate module $\text{Pl}\left(\Delta_r^n\right)$ are isomorphic as $S_n$-modules.
\end{thm}

\begin{defn}
	Let $\mathcal{I}=\mathcal{I}_r^{n}$ denote the set of all $(i_1,\ldots, i_n)\in\{0,\ldots, r-1\}^n$ with $ \sum_{j=0}^{r-1}i_j\equiv 1 \mod r$.
	
	For each $n$-tuple $I=(i_1,\ldots, i_n)\in\mathcal{I}$, let
	$$\epsilon_{I}=\frac{1}{r^n}\prod_{j=1}^n\left(\sum_{k=0}^{r-1} \left(q^{-i_j}e_j\right)^k\right).$$
\end{defn}

Remark that the condition $\sum_k i_k\equiv 1$ is dual to $e_1\cdots e_n=q$, and when both conditions are applied, the remaining variables are changed by a Fourier transform.

\begin{thm}\label{IdempotentsDiagonalize}
	The elements $\epsilon_I$ form a partition of unity of $\mathcal{C}_r^n$ into $r^{n-1}$ one-dimensional subspaces, that is, 
	\begin{eqnarray*}
		\epsilon_I \epsilon_{I'} & = & \delta_{I,I'}\epsilon_I\\
		\sum_{I\in\mathcal{I}} \epsilon_I  & = & \mathbf{1}.
	\end{eqnarray*}
	These subspaces are simultaneous eigenvectors for multiplication by $e_i$, that is, 
	$$e_j \epsilon_I=q^{i_j}\epsilon_I,$$
	and so 
	$$(e_1\cdots e_n)\epsilon_I=q^{i_1+\cdots+i_n}\epsilon_I=q\epsilon_I.$$
	
\end{thm}
The character formula for the action of the symmetric group on $\mathcal{C}_r^n$, and by Theorem \ref{permutohedronAlgIsomorphism} also on $\text{Pl}(\Delta_r^n)$, can be obtained by solving a \textit{modular} Diophantine equation, as follows.

\begin{prop}
	The trace of a permutation of cycle type $\lambda_1,\ldots, \lambda_k$ acting on $\mathcal{C}_r^n$ is equal to the number of solutions to the equation 
	$$\sum_{i=1}^k \lambda_i x_i\equiv 1\text{ mod $r$}.$$
This number is equal to $r^{k-1}$ if $\gcd(\lambda_1,\ldots, \lambda_k,r)=1$ and is 0 otherwise.  
\end{prop}

\begin{example}
	In Figure \ref{fig:idempotent-character} the fixed points of the cycles $(12)$ and $(123)$ acting on the set $\mathcal{I}_{10}^3$ are indicated with colored dots.  The light green dot is fixed by $(123)$, while the green together with the nine red dots are fixed by $(12)$.  The value of the character evaluated on a permutation $\sigma$ is thus given by the number of fixed points of $\sigma$ on this set.  
\begin{figure}[h!]
	\centering
	\includegraphics[width=.35\linewidth]{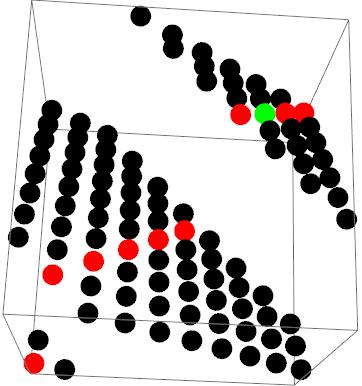}
	\caption{}
	\label{fig:idempotent-character}
\end{figure}

\end{example}

\subsection{Multiplicities of the Irreducible Representations}
It is an obvious question to ask for the multiplicities of the irreducible representations in $\text{Pl}(\Delta_r^n)$ and $\text{Pl}(B_{a,b})$.  We have computed the multiplicities for the trivial representation in $\text{Pl}(\Delta_r^n)$, for which the multiplicities are enumerated by Lyndon words, but a full combinatorial and geometric interpretation of the general case, for arbitrary irreducible representations, is beyond the scope of this note and shall appear elsewhere.   

The first few cases of the multiplicities of the irreducible representations in $\text{Pl}(\Delta_r^n)$ are given in the table below.  Here for example $4\cdot 2^1 1^2$ means that the irreducible $S_4$ representation labeled by the partition $2+1+1$ of $4$ occurs with multiplicity $4$.

$$\begin{tabular}{|c|c|c|c|c|}
\hline 
$n\setminus r $& 1 & 2 & 3 & 4 \\ 
\hline 
1& $1^1$ &$ 1^1$& $1^1 $& $1^1$ \\ 
\hline 
2& $2^1$ & $2^1+1^2$ & $2\cdot 2^1+1^2$ & $2\cdot 2^1+2\cdot 1^2$\\ 
\hline 
3&$ 3^1$ & $2\cdot 3^1+2^1 1^1$ & $3\cdot 3^1+3\cdot 2^1 1^2$ & $5\cdot 3^1+5\cdot 2^1 1^1+1^3$  \\ 
\hline 
4& $4^1$ & $2\cdot 4^1+2\cdot 3^1 1^1$ & $5\cdot 4^1+2\cdot 2^2+5\cdot 3^1 1^1+2^1 1^2$ &$8\cdot 4^1+4\cdot 2^2+12\cdot 3^1 1^1+4\cdot 2^1 1^2$  \\ 
\hline 
\end{tabular} $$

By explicit enumeration of fixed elements in the idempotent basis, we found the generating functions for the multiplicities of the trivial representations to be
$$\frac{1}{(1-x)^3 \left(1+x+x^2\right)}=1+2 x+3 x^2+5 x^3+7 x^4+9 x^5+\cdots $$
$$\frac{1}{(1-x)^4 (1+x)^2}=1+2 x+5 x^2+8 x^3+14 x^4+20 x^5+\cdots,$$
where the coefficients of the series can be seen in columns labeled $3$ and $4$ in the table.

\section{Acknowledgements}
I am grateful to Adrian Ocneanu for sharing his conjectures and encouraging me to develop the theory as much as possible, and for countless hours of intensive discussions.  I thank Sergei Tabachnikov and Nigel Higson for comments and suggestions for the final draft of my thesis.  I thank Alexander Postnikov and Sinai Robins for stimulating conversations.

\end{document}